\documentclass[11pt]{amsart}  
\usepackage{amsmath,amssymb}
\usepackage[dvips]{graphicx}  
\usepackage{psfrag}
\def\COMMENT#1{} 
\def\TASK#1{}

\def\noproof{{\unskip\nobreak\hfill\penalty50\hskip2em\hbox{}\nobreak\hfill%
        $\square$\parfillskip=0pt\finalhyphendemerits=0\par}\goodbreak}
\def\endproof{\noproof\bigskip}
\newdimen\margin   
\def\textno#1&#2\par{%
    \margin=\hsize
    \advance\margin by -4\parindent
           \setbox1=\hbox{\sl#1}%
    \ifdim\wd1 < \margin
       $$\box1\eqno#2$$%
    \else
       \bigbreak
       \hbox to \hsize{\indent$\vcenter{\advance\hsize by -3\parindent
       \sl\noindent#1}\hfil#2$}%
       \bigbreak
    \fi}
\def\proof{\removelastskip\penalty55\medskip\noindent{\bf Proof. }}

\def\eps{\varepsilon}

\newtheorem{firstthm}{Proposition}

\newtheorem{prop}[firstthm]{Proposition}

\newtheorem{conj}[firstthm]{Conjecture}

\newtheorem{question}[firstthm]{Question}

\begin{document}
\title{A note on some embedding problems for oriented graphs}
\author{Andrew Treglown}
\thanks {The author was supported by the EPSRC, grant no.~EP/F008406/1.}
\date{\today}

\begin{abstract}
We conjecture that every oriented graph $G$ on $n$ vertices with $\delta ^+ (G) , \delta ^- (G) \geq 5n/12$ contains the square of a Hamilton
cycle. We also give a conjectural bound on the minimum semidegree which ensures a perfect packing of transitive 
triangles in an oriented graph. A link between Ramsey numbers and perfect packings of transitive tournaments is also considered.
\end{abstract}
\maketitle
\subsection{Powers of Hamilton cycles} 
One of the most studied problems in graph theory concerns finding sufficient conditions that ensure a graph contains a Hamilton cycle. 
Dirac~\cite{dirac} showed that any graph $G$ on $n \geq 3 $ vertices has a Hamilton cycle provided that it has minimum degree $\delta (G)$ at least 
$n/2$. For a digraph $G$ it is natural to consider its \emph{minimum semidegree} $\delta ^0 (G)$, which 
is the minimum of its minimum outdegree $\delta ^+ (G)$ and its minimum indegree $\delta ^- (G)$.
(The digraphs we consider do not have loops and we allow at most one edge in each
direction between any pair of vertices.)
Ghouila-Houri~\cite{gh} proved that every digraph
$G$ on $n\geq 2$ vertices with $\delta ^0 (G)\geq n/2$ is Hamiltonian. 

An important subclass of digraphs is the class of
 oriented graphs: these are the digraphs which do not contain any 2-cycles. Keevash, K\"uhn and Osthus~\cite{kko}
showed that any sufficiently large oriented graph $G$ on $n$ vertices with $\delta ^0 (G) \geq (3n-4)/8$ is Hamiltonian, thereby proving a conjecture 
of H\"aggkvist~\cite{HaggkvistHamilton}. For a detailed account of other such results concerning 
Hamilton cycles in directed and oriented graphs see~\cite{cyclesurvey}.

A generalisation of the notion of a Hamilton cycle is that of the $r$th power of a Hamilton cycle. Indeed, the $r$th power of a Hamilton cycle
$C$ is obtained from $C$ by adding an edge between every pair of vertices of distance at most $r$ on $C$. Seymour~\cite{seymour} conjectured the
following strengthening of Dirac's theorem.
\begin{conj}[Seymour~\cite{seymour}]\label{seymour}
Let $G$ be a graph on $n$ vertices. If $\delta (G) \geq \frac{r}{r+1} n$
then $G$ contains the $r$th power of a Hamilton cycle.
\end{conj}
P\'osa~(see \cite{posa}) had earlier proposed the conjecture in the case of the square of a Hamilton cycle (that is, when $r=2$). 
Koml\'os, S\'ark\"ozy and Szemer\'edi~\cite{kss} proved Conjecture~\ref{seymour} for sufficiently large graphs.

The notion of the $r$th power of a Hamilton cycle also makes sense in the digraph setting: 
In this case the $r$th power of a Hamilton cycle $C$ is the digraph obtained from
$C$ by adding a directed edge from $x$ to $y$ if there is a path of length at most $r$ from $x$ to $y$ on $C$.
Bollob\'as and H\"aggkvist~\cite{bolhag} proved that given any $\eps>0$ and any $r \in \mathbb N$, all sufficiently large tournaments 
$T$ on $n$ vertices with $\delta ^0 (T) \geq (1/4+\eps)n$ contain the $r$th power of a Hamilton cycle.

One would expect that the minimum semidegree threshold that ensures a digraph contains the $r$th power of a Hamilton cycle is the `same' as the
condition in Conjecture~\ref{seymour}. But it is far less clear at first sight what to expect in the oriented case. We propose the following
oriented graph analogue of P\'osa's conjecture.
\begin{conj}\label{co2} Suppose $G$ is an oriented graph on $n$ vertices such that $\delta ^0 (G) \geq 5n/12.$ Then $G$ contains the square
of a Hamilton cycle.
\end{conj}
The following proposition shows that, if true, Conjecture~\ref{co2} is `best possible'. 
\begin{prop}\label{prop1} Let $n \in \mathbb N$ be divisible by $12$. Then there is an oriented graph $G$ on $n$ vertices with $\delta ^0 (G)= 5n/12 -1$ which
does not contain the square of a Hamilton cycle.
\end{prop}
\proof Let $G$ denote the oriented graph on $n$ vertices whose vertex set consists of the sets $A,B,C,D$ and $E$ where $|A|=n/6+1$, $|B|=n/6-1$,
$|C|=n/3$ and $|D|=|E|=n/6$. The edge set of $G$ is obtained as follows: Add all possible edges from $A\cup B$ to $C$, from $C$ to $D\cup E$, from
$D$ to $A\cup B$ and from $E$ to $A \cup D$. Let $B$, $C$ and $D$ all induce tournaments that are as regular as possible (so $\delta ^0 (G[B])=
\delta ^0 (G[D])= n/12-1$ and $\delta ^0 (G[C])=n/6-1$). We add edges between $A$ and $B$ in such a way that every vertex in $A$ sends and receives
at least $n/12-1$ edges to and from $B$, and every vertex in $B$  sends and receives
at least $n/12$ edges to and from $A$. Similarly, we add edges between $B$ and $E$ in such a way that every vertex in $B$ sends and receives
 $n/12$ edges to and from $E$, and every vertex in $E$  sends and receives
at least $n/12-1$ edges to and from $B$. $A$ and $E$ are both independent sets (see Figure~1). So $\delta ^0 (G) = 5n/12-1$. 

Assume that $G$ contains the square of a Hamilton cycle $F$. Since $|B|<|E|$, showing that $F$ 
must visit $B$ between any two visits of $E$ would yield a contradiction. Thus, consider any vertex $e\in E$. Its predecessor $c_1$ on $F$ lies in 
$B \cup C$, so without loss of generality we may assume that $c_1 \in C$. The predecessor $c_2$ of $c_1$ on $F$ must lie in
$N^- (e)\cap N^- (c_1) \subseteq B \cup C$. So without loss of generality we may assume that $c_2 \in C$. The predecessor $c_3$ of $c_2$ on $F$
lies in $A \cup B \cup C$. Again we are done if $c_3 \in B$, so we assume that $c_3 \in A \cup C$. Since $F$ visits all the vertices of $G$ we must 
eventually arrive at a predecessor $a \in A$ whose successor $c$ on $F$ lies in $C$. But now the predecessor of $a$ on $F$ must lie
in $N^- (c) \cap N^- (a) \subseteq B$, as required.
\endproof
\begin{figure}[htb!] 
\begin{center}\footnotesize 
\psfrag{1}[][]{$B$}
\psfrag{2}[][]{$D$}
\psfrag{3}[][]{$A$}
\psfrag{4}[][]{$E$} 
\psfrag{5}[][]{$C$}
\includegraphics[width=0.65\columnwidth]{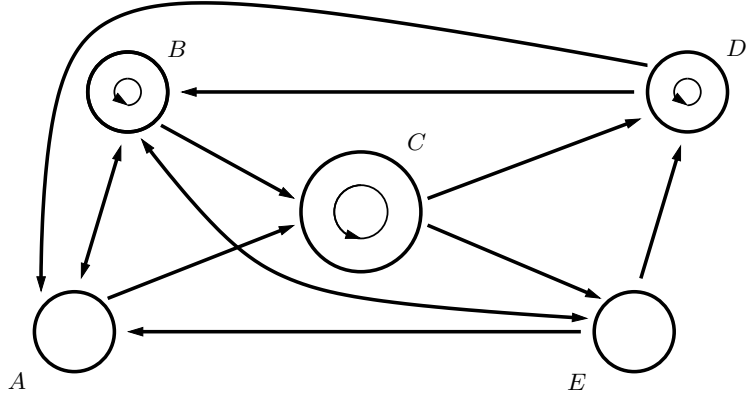}  
\caption{The oriented graph $G$ from Proposition~\ref{prop1}} 
\end{center}
\end{figure}
\subsection{Transitive triangle packings} Given an (oriented) graph $H$, a \emph{perfect $H$-packing} in an (oriented) graph $G$ is a collection of
vertex-disjoint copies of $H$ which covers all the vertices in $G$. (Perfect $H$-packings are also referred to as \emph{$H$-factors} or \emph{perfect
$H$-tilings}.)
Perfect $H$-packings in graphs have been widely studied. Hajnal and Szemer\'edi~\cite{hs} characterised the minimum degree which ensures a graph $G$
contains a perfect $K_r$-packing. More recently, K\"uhn and Osthus~\cite{kuhn2} characterised, up to an additive constant, the minimum degree which
ensures a graph $G$ contains a perfect $H$-packing for an arbitrary graph $H$. Far less is known in the oriented graph case.
Keevash and Sudakov~\cite{keevs} 
showed that any oriented graph $G$ on $n$ vertices with $\delta ^0 (G) \geq(1/2-o(1))n$ contains a packing of cyclic triangles covering all but at most $3$ vertices. 

It is natural to ask for the minimum semidegree of an oriented graph which ensures a perfect packing of transitive triangles $T_3$. Note that
if $3$ divides $|G|$ then a necessary condition for an
oriented graph $G$ to contain a square of a Hamilton cycle is that $G$ contains a perfect packing of transitive triangles. 
Let $\delta (G)$ denote the minimum degree of an oriented graph $G$ (that is, the minimum number of edges incident to a vertex in $G$). 
The following 
proposition from~\cite{yuster} implies that a minimum semidegree as in Conjecture~\ref{co2} ensures a perfect $T_3$-packing.
\begin{prop}[Yuster~\cite{yuster}]\label{prop2} Suppose $G$ is an oriented graph whose order $n$ is divisible by $3$. If $\delta (G)
\geq 5n/6$ then $G$ contains a perfect $T_3$-packing.
\end{prop} 
Proposition~\ref{prop2} is best possible in the sense that there are oriented graphs $G$ whose order $n$ is divisible by $3$ and where $\delta (G)=
(5n-3)/6$ 
but which do not contain a perfect $T_3$-packing. (Indeed, consider the oriented graph $G$ on $6m+3$ vertices consisting of $3$ vertex sets $A$, $B$
and $C$ where $|A|=|B|=m+1$ and $|C|=4m+1$, and such that $C$ induces a tournament, $A$ sends out all possible edges to $B$, $B$ sends out all possible
edges to $C$ and $C$ sends out all possible edges to $A$. Then $G$ does not contain a perfect $T_3$-packing since every copy of $T_3$ in $G$ has at most one vertex in $A \cup B$.) However, when considering embeddings in oriented graphs, it seems that the more natural parameter to look at
is the minimum semidegree.
We believe that, in terms of minimum semidegree, one can improve on the bound given in Proposition~\ref{prop2}.
\begin{conj}\label{co3} Suppose $G$ is an oriented graph whose order $n$ is divisible by $3$. If $\delta ^0 (G) \geq 7n/18$ then $G$ contains
a perfect $T_3$-packing.
\end{conj}
If true, Conjecture~\ref{co3} would characterise the minimum semidegree which ensures an oriented graph has a perfect $T_3$-packing.
\begin{prop}\label{triangle} Let $n\in \mathbb N$ be divisible by $18$. Then there is an oriented graph $G$ on $n$ vertices with
$\delta ^0 (G) = 7n/18-1$ which does not contain a perfect $T_3$-packing.
\end{prop}
\proof Let $G$ denote the oriented graph on $n$ vertices whose vertex set consists of the
sets $A$, $B$, $C$ and $D$ where $|A|=2n/9 +1$, $|B|=|C|=2n/9$ and $|D|=n/3-1$ and whose edge set is obtained as follows: Add all 
possible edges from $A$ to $B$, from $B$ to $C$ and from $C$ to $A$. 
Let $D$ induce a regular tournament.
Partition $D$ into two sets $D'$ and $D''$ of sizes $n/6$ and $n/6-1$
respectively. Add all possible edges from $D'$ to $B \cup C$, from $A$ to $D'$, from $D''$ to $A$ and from $B \cup C$ to $D''$ (see Figure~2). 
It is 
easy to see that $\delta ^0 (G) =7n/18 -1$. Note that $G$ does not have a perfect $T_3$-packing since 
every copy of $T_3$ in $G$ must have at least one vertex in $D$.
\endproof 
\begin{figure}[htb!] 
\begin{center}\footnotesize  
\psfrag{1}[][]{$B$}
\psfrag{2}[][]{$C$}
\psfrag{6}[][]{$A$}
\psfrag{7}[][]{$D$}
\includegraphics[width=0.35\columnwidth]{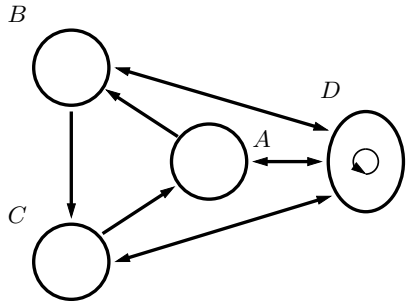}\end{center}
\caption{The oriented graph $G$ from Proposition~\ref{triangle}} 
\end{figure}
\subsection{Packing transitive tournaments} Let $T_k$ denote the transitive tournament on $k$ vertices. In light of Conjecture~\ref{co3} we ask
the following question.
\begin{question}\label{q1} What minimum semidegree condition ensures that an oriented graph contains a perfect $T_k$-packing?
\end{question}
Recall that in the oriented graph $G$ given in Proposition~\ref{triangle} the vertex set $A \cup B \cup C$ induces an oriented graph which does not
contain a copy of $T_3$. This is the `reason' why $G$ does not contain a perfect $T_3$-packing. It would be of interest to establish whether
the extremal examples, in terms of perfect $T_k$-packings, take a similar form. Thus, Question~\ref{q1} is closely linked to the following
question.
\begin{question}\label{q2} What minimum semidegree condition ensures that an oriented graph contains a copy of $T_k$?
\end{question}
Valadkhan~\cite{val} has investigated this problem with respect to density conditions (a wider collection of problems of this nature are considered in~\cite{brown}).
It is easy to see that an oriented graph $G$ on $n$ vertices with $\delta ^0 (G) >n/3$ contains a copy of $T_3$ (and the blow-up
of a cyclic triangle shows that this bound is best 
possible). 
\subsection{Perfect packings and Ramsey numbers}
The \emph{oriented tiling Ramsey number $\overrightarrow{TR}(k)$ of $k$} is the smallest integer $n$ divisible by $k$ such that
any orientation of the complete graph $K_n$ contains a perfect $T_k$-packing. 
Erd\H{o}s (see~\cite{reid})
proved the existence of these numbers.
The following simple result gives a bound on the minimum
degree which ensures an oriented graph $G$ contains a perfect $T_k$-packing.
\begin{prop}\label{yusterprop} Suppose $G$ is an oriented graph whose order $n$ is divisible by $k$ and such that $\delta (G) \geq 
(1-\frac{1}{\overrightarrow{TR}(k)})n.$ Then $G$ contains a perfect $T_k$-packing.\COMMENT{don't need +constant here!}
\end{prop}
\noindent
{\bf Sketch proof.} Let $m:={\overrightarrow{TR}(k)}$.
Consider the case when $m$ divides $n$.
By disregarding the orientations of the edges of $G$ we obtain
a graph $G^*$ on $n$ vertices with $\delta (G^*) \geq (1-\frac{1}{m})n$. The Hajnal-Szemer\'edi theorem~\cite{hs} implies that $G^*$ has a perfect $K_m$-packing. By definition of $m$ this implies that $G$ has a perfect $T_k$-packing. If $n$ is not divisible by $m$, we remove
a number of vertex-disjoint copies of $T_k$ from $G$ until $m$ divides $|G|$. We then proceed as before.
\endproof
Note that $\overrightarrow{TR}(3)=6$ so Proposition~\ref{yusterprop} implies Proposition~\ref{prop2}.
In view of Proposition~\ref{yusterprop} it is natural to seek good upper bounds on $\overrightarrow{TR}(k)$. The \emph{oriented Ramsey number
$\overrightarrow{R}(k)$ of $k$} is the smallest integer $n$ such that any orientation of $K_n$ contains a copy of $T_k$. 
The following proposition gives an upper bound on $\overrightarrow{TR}(k)$ in terms of oriented Ramsey numbers.
\begin{prop} Given any $k \in \mathbb N$, $\overrightarrow{TR}(k) \leq \overrightarrow{R}(2k-1)+(2k-1)\overrightarrow{R}(k)$.
\end{prop}
\proof We use the same trick as Caro used in~\cite{caro}. Let $n$ be the largest integer divisible by $k$ such that $n \leq 
\overrightarrow{R}(2k-1)+(2k-1)\overrightarrow{R}(k)$ and $\ell$ the largest integer divisible by $k$ which satisfies $\ell \leq 
\overrightarrow{R}(k)$. Consider any orientation $\overrightarrow{K}$ of $K_n$. By definition of $n$, $\overrightarrow{K}$ contains
$\ell$ vertex-disjoint copies of $T_{2k-1}$. We can cover all but $\ell$ of the remaining vertices of $\overrightarrow{K}$ with
vertex-disjoint copies of $T_k$. Each of the $\ell$ uncovered vertices $x$ are paired off with one of our copies $T'_{2k-1}$
of $T_{2k-1}$. Since $x$ either sends out at least $k$ edges to $T'_{2k-1}$ in $\overrightarrow{K}$  or receives at least $k$ edges from
$T'_{2k-1}$ in $\overrightarrow{K}$, we have that the oriented subgraph of $\overrightarrow{K}$ induced by $V(T'_{2k-1})\cup \{x\}$ contains a perfect
$T_k$-packing. Thus $\overrightarrow{K}$ contains a perfect $T_k$-packing.
\endproof 
The numbers $\overrightarrow{R}(k)$ are known for $k \leq 6$ (see~\cite{reidparker,san1}).
Sanchez-Flores~\cite{san} showed that $\overrightarrow{R}(7)\leq 54$ which by an induction argument implies that 
$\overrightarrow{R}(k)\leq 54 \cdotp 2^{k-7}$ for $k\geq 7$ (this is the best known
general upper bound on oriented Ramsey numbers). Note also that $\overrightarrow{R}(k)\leq R(k)$ where $R(k)$ denotes the Ramsey number of $k$.

\section*{Acknowledgement}
The author would like to thank Daniela K\"uhn, Richard Mycroft and Deryk Osthus for helpful comments and discussions.

\medskip 

{\footnotesize \obeylines \parindent=0pt

Andrew Treglown
School of Mathematics
University of Birmingham
Edgbaston
Birmingham
B15 2TT
UK
}
\begin{flushleft}
{
\tt treglowa@maths.bham.ac.uk}
\end{flushleft}\end{document}